\documentclass{amsart}
\usepackage[latin1]{inputenc}
\usepackage{amssymb}
\usepackage{amsmath}
\usepackage{enumerate}
\usepackage{epsfig,color}
\usepackage[all]{xy}

\newtheorem{theorem}{Theorem}[section]
\newtheorem{lemma}[theorem]{Lemma}
\newtheorem{proposition}[theorem]{Proposition}
\newtheorem{corollary}[theorem]{Corollary}

\newtheorem{definition}[theorem]{Definition}

\newtheorem{remark}[theorem]{Remark}
\newtheorem{rem}[theorem]{Remark}

\numberwithin{equation}{section}
\newcommand{\rk}{\mbox{rank}}

\newcommand{\MW}{\mbox{MW}}

\newcommand{\Aut}{\mbox{Aut}}
\newcommand{\id}{\mbox{id}}

\newcommand{\ra}{\rightarrow}

\newcommand{\Z}{\mathbb{Z}}

\newcommand{\contr }{ \mbox{contr}}
\newcommand{\Cn}{C_{(n)}}
\newcommand{\Yns}{Y_{(n)}^{(2)}}
\newcommand{\Xns}{X_{(n)}^{(2)}}
\newcommand{\Ynt}{Y_{(n)}^{(3)}}
\newcommand{\Xnt}{X_{(n)}^{(3)}}
\newcommand{\Ynl}{Y_{(n)}^{(l)}}
\newcommand{\Znl}{Z_{(n)}^{(l)}}
\newcommand{\Xnl}{X_{(n)}^{(l)}}
\newcommand{\Gnl}{G_{(n)}^{(l)}}

\makeatletter
\def\blfootnote{\xdef\@thefnmark{}\@footnotetext}

\title[Elliptic fibrations and Kodaira dimensions of Schreieder's varieties]{Elliptic fibrations and Kodaira dimensions of Schreieder's varieties}

\author{Alice Garbagnati}

\address{Alice Garbagnati, Universit\`a degli Studi di Milano, Dipartimento di Matematica, via Cesare Saldini 50 20133 Milano, Italy }
\email{alice.garbagnati@unimi.it}
\urladdr{https://sites.google.com/site/alicegarbagnati/}

\begin{document}
	
	\begin{abstract}
		We discuss the birational geometry and the Kodaira dimension of certain varieties previously constructed by Schreieder, proving that in any dimension they admit an elliptic fibration and they are not of general type. The $l$-dimensional variety $\Ynl$, which is the quotient of the product of a certain curve $\Cn$ by itself $l$ times by a group $G\simeq\left(\Z/n\Z\right)^{l-1}$ of automorphisms, was constructed by Schreieder to obtain varieties with prescribed Hodge numbers. If $n=3^c$ Schreieder constructed an explicit smooth birational model of it, and Flapan proved: the Kodaira dimension of this smooth model is 1, if $c>1$; if $l=2$ it is a modular elliptic surface; if $l=3$ it admits a fibration in K3 surfaces. 
		
		In this paper we generalize these results: without any assumption on $n$ and $l$ we prove that $\Ynl$ admits many elliptic fibrations and its Kodaira dimension is at most 1. Moreover, if $l=2$ its minimal resolution is a modular elliptic surface, obtained by a base change of order $n$ on a specific extremal rational elliptic surface; if $l\geq 3$ it has a birational model which admits a fibration in K3 surfaces and a fibration in $(l-1)$-dimensional varieties of Kodaira dimension at most 0. 
	\end{abstract}
	\maketitle
	
	\section{Introduction}\label{intro}
	
	The aim of this paper is to reconsider certain varieties constructed by Schreieder as quotient of the product of a very special curve, $\Cn$, with itself $l$-times.
	The curve $\Cn$ is a hyperelliptic curve $\Cn\ra\mathbb{P}^1_u$ endowed with a special symmetry among the ramification points. Denoted by $u$ an affine coordinate of $\mathbb{P}^1$, the equation of $\Cn$  as double cover of $\mathbb{P}^1$, is $v^2=u^n-1$. So the map $(v,u)\ra (v, \zeta_n u)$, where $\zeta_n$ is an $n$-primitive root of unity, induces an automorphism of $\Cn$, denoted by $\alpha$. 

	For $n=3,4$, $\Cn$ is the elliptic curve with complex multiplication of order $n$.
	
	For each $l$, and $n=3,4$ Cynk and Hulek proved that a quotient of the product of $\Cn$ by itself $l$-times by a specific group of automorphisms admits a resolution which is a modular Calabi--Yau $l$-fold, see \cite{CH}.
	
	More recently a generalization of the construction by Cynk and Hulek was considered by Schreieder (see \cite{S}): let $\Cn^{\times l}$ be the product of $\Cn$ by itself $l$-times and $G\simeq \left(\Z/n\Z\right)^{l-1}$ be a specific group of automorphisms of $\Cn^{\times l}$, then the variety $\Ynl$ is defined as the (singular) quotient  $\Cn^{\times l}/G$. 
	
	For the varieties $\Ynl$ such that $n=3^c$, Schreieder also constructed explicit resolutions: each variety $Y_{(3^c)}^{(l)}$ is a quotient of $Y_{(3^c)}^{(l-h)}\times Y_{(3^c)}^{(h)}$, $h<l$, and so one obtains a desingularization which generalizes the construction in \cite{CH}. In \cite{F} Flapan proved that the smooth varieties constructed by Schreieder have Kodaira dimension 1 if $c>1$, the 2-dimensional ones are elliptic modular surfaces and the 3-dimensional ones admit a fibration in K3 surfaces.
	
	We give generalizations of the Flapan's results without any assumption on $n$ and on $l$, indeed we obtain elliptic fibrations (and K3-fibrations if $l>2$) defined on birational models of all the varieties $\Ynl$, and we prove  $\Ynl$ have Kodaira dimension, $k(\Ynl)$, not bigger than 1. In dimension 2 and 3 we fully recover the statements of \cite{F} without the assumption on $n$. These results can be achieved without constructing an explicit desingularization of $\Ynl$, but just by considering the equations of certain (possiby singular) birational models.  
	
	To state the main theorem one has to fix a definition of  ``elliptic fibration" also for varieties with dimension higher than 2. We will say that a surjective map $\mathcal{E}:X\ra B$ from a possibly singular $l$-dimensional variety $X$ to a smooth $(l-1)$-dimensional variety $B$ is an {\it elliptic fibration} if $\mathcal{E}$ has connected fibers, its generic fiber is a smooth genus 1 curve and there is a rational map $s:B\dashrightarrow X$ such that $\mathcal{E}\circ s$ is the identity where it is defined. So $X$ admits an elliptic fibration if it is birational to a variety defined by a Weierstrass equation, over $B$.
	
	The synthesis of the results in the paper is the following:\\
	
	{\bf Main Theorem }{\it Let $n$ and $l$ be two integers bigger than 1. Then $\Ynl$ is birational to a variety $\Znl$ which admits an elliptic fibration with basis $(\mathbb{P}^1)^{\times (l-1)}$.
		\begin{itemize}\item If $l=2$ and $n>2$ (resp. $n=2$) the minimal model of $\Ynl$ (resp. a smooth  birational model of $Y_{(2)}^{(l)}$) is a modular elliptic surface. It is a rational surface if $n=2$, a K3 surface if $n=3,4$, a properly elliptic surface (i.e. $k\left(\Yns\right)=1$) if $n>4$.
			\item If $l
			\geq 3$, $\Ynl$ admits a fibration in K3 surfaces with generic Picard number 19.  For every $n$, $k\left(\Ynl\right)\leq 1$ and in particular if $n=2$, $k\left(\Ynl\right)=-\infty$, if $n=3,4$, $k\left(\Ynl\right)=0$.
	\end{itemize}}
	
	To prove the main theorem one first observes that the varieties $\Ynl$ are birational to a $\Z/n\Z$ quotient of $Y_{(n)}^{(l-1)}\times \Cn$. This naturally gives an isotrivial fibration $\Ynl\ra\mathbb{P}^1$ whose fibers are isomorphic to $Y_{(n)}^{(l-1)}$, thus the knowledge of properties (and in particular of a nice equation) of $Y_{(n)}^{(l-1)}$ gives information also on $\Ynl$. The definition of the varieties $\Ynl$ and the analysis of the relations between these varieties for different values of $n$ and $l$ are the topics of Section \ref{sec: Peeliminaries}.
	
	In Section \ref{sec: Yns} we consider the first interesting case, i.e. the 2-dimensional variety $\Yns$. In dimension 2, one is able to prove that a smooth model of $\Yns$ (indeed the minimal model for each $n>2$) is an elliptic modular surface, and in particular it is not of general type. The elliptic fibration is completely described (both the reducible fibers and the Mordell--Weil group) and the full Hodge diamond of the minimal model is computed, see Theorem \ref{theorem surfaces} and Proposition \ref{prop: the surface Z_n}. These results are achieved in two different ways: either one constructs the minimal resolution of $\Yns$ and finds explicitly a nef  divisor whose associated map is the required elliptic fibration (see Subsection \ref{subsec: geometry of minimal resolution surface}), or one considers a very particular equation of the quotient surface $\Yns=\left(\Cn\times \Cn\right)/(\Z/n\Z)$ (see Proposition \ref{prop: Yns birat to base change}) which exhibits directly the elliptic fibration. In both the contexts one observes that $\Yns$ is an $n:1$ cover of $\mathbb{P}^1\times\mathbb{P}^1$ branched on two (singular) curves of bidegree $(2,2)$ and indeed it is birational to a $n:1$ cover of a rational elliptic surface, $R$ (which does not depend on $n$).
	
	In Section \ref{sec: 3-fold} we consider the 3-dimensional variety $\Ynt$. We give an equation of $\Ynt$ (obtained from the one of $\Yns$), which shows that $\Ynt$ admits two different fibrations: one is the elliptic fibration with basis $\mathbb{P}^1\times\mathbb{P}^1$ mentioned in the main theorem; the other is a fibration in K3 surfaces with basis $\mathbb{P}^1$. Hence also in the 3-dimensional case, $\Ynt$ admits a fibration over $\mathbb{P}^1$ in codimension 1 subvarieties with trivial canonical bundle, in particular the Kodaira dimension of $\Ynt$ is not higher than 1.
	
	In Section \ref{sec: l-fold} we iterate the process applied in the previous sections to find an equation for $\Ynl$, and in particular of fibrations on it. One finds the elliptic fibration and the K3-fibration mentioned in the main theorem and a fibration in varieties of codimension 1, whose generic fibers do not depend on $n$ and have Kodiara dimension at most 0. This latter fibration is the analogous of the elliptic fibration on $\Yns$ and of the K3-fibrations on $\Ynt$ and is used to bound the Kodaira dimension of $\Ynl$.

\subsection*{Acknowledgements} 
I'm grateful to Bert van Geemen for his enlightening comments. The author would like to thank Matthias Sch\"utt  and Luca Tasin and for useful conversations.	
\section{Preliminaries: the curves $\Cn$ and the varieties $\Ynl$}\label{sec: Peeliminaries}
In this section we define the varieties  $\Ynl$  and we present some of their properties. These are varieties of dimension $l$, where $l$ is arbitrarily big, obtained by the quotient of the product of a certain curve $\Cn$ (defined in Section \ref{subsec: The curve Cn}) $l$ times by itself. The construction of $\Ynl$ is presented in Section \ref{subsect: The varieties Ynl}. In Section \ref{subsec: relations between Ynl Ymh} we describe the relations between $Y_{(n)}^{(l)}$ and $Y_{(m)}^{(h)}$, for $n\neq m$ and/or $l\neq h$ and, in particular, the ones between the Kodaira dimension of the varieties $\Ynl$ for different values of $l$ (in Proposition \ref{prop: Kodaira dimension}).

\subsection{The curves $\Cn$}\label{subsec: The curve Cn} The curve $\Cn$ is obtained both as $2:1$ cover and as $n:1$ cover of $\mathbb{P}^1$. We chose to define it as a very special $n:1$ cover of $\mathbb{P}^1$. 
\begin{definition}
	
	If $n$ is odd the curve $C_{(n)}$ is the $n:1$ cover of $\mathbb{P}^1$ totally branched over 3 points such that $C_{(n)}$ admits an involution $\iota_n$ switching two ramification points and fixing the third.
	
	If $n$ is even the curve $C_{(n)}$ is the $n:1$ cover of $\mathbb{P}^1$ totally branched on 2 points and branched with multiplicity $n/2$ on a third point such that $C_{(n)}$ admits an involution $\iota_n$ switching the two points of total ramification and the two points of partial ramification.
\end{definition}

An equation of $C_{(n)}$ as $n:1$ cover of $\mathbb{P}^1_{(v:w)}$ is \begin{equation}\label{eq: equation Cn}u^n=(v^2-w^2)w^{n-2},\end{equation}
with cover map $f_{(n)}:\Cn\ra\mathbb{P}^1_{(v:w)}$. 

The involution $\iota_n:\Cn\ra \Cn$ is $\iota_n:(u,(v:w))\ra(u,(-v:w))$ and the cover automorphism  
 \begin{equation}\label{eq: equation alphan}\alpha_n:\Cn\ra \Cn\mbox{ is }\alpha_n:(u,(v:w))\ra(\zeta_nu,(v:w))\end{equation}
where $\zeta_n$ is a primitive $n$-th root of unity. We assume that, if $m|n$, then $\zeta_m=\zeta_n^{n/m}$.

By considering the affine equation of $\Cn$ obtained putting $w=1$, one finds the equation $u^n=v^2-1$, which exhibits $\Cn$ as the hyperelliptic curve $v^2=u^n+1$, $2:1$ cover of $\mathbb{P}^1_u$. Then $\iota_n$ is the hyperelliptic involution. 

Here we collect some very well known results on the curve $\Cn$ and on its automorphism which will be used in the following. If $n$ is odd, we denote by  $P_1$, $P_{-1}$, $P_{\infty}$ the three points of $\Cn$ fixed by $\alpha_n$, which are respectively mapped to $(1:1)$, $(-1:1)$, $(1:0)$ by $f_{(n)}:\Cn\ra\mathbb{P}^1_{(v:w)}$.

If $n$ is even, we denote by  $P_1$ and $P_{-1}$ the two points fixed by $\alpha_n$, respectively mapped to $(1:1)$, $(-1:1)$ by $f_{(n)}:\Cn\ra\mathbb{P}^1_{(v:w)}$, and by $P_{\infty}^1$ and $P_{\infty}^2$ the two points switched by $\alpha_n$ and mapped to $(1:0)$ by $f_{(n)}:\Cn\ra\mathbb{P}^1_{(v:w)}$.

\begin{proposition}\label{prop: properties curve Cn}
The curve $\Cn$ has genus 
$$g(\Cn)
=\left\{\begin{array}{llll}\frac{n-1}{2}&\mbox{ if }n\mbox{ is odd},\\\frac{n}{2}&\mbox{ if }n\mbox{ is even.}
\end{array}\right.$$ A basis of $H^{1,0}(\Cn)$ is $\{\omega_i:=u^{i-1}du/v,\ i=1,\ldots, g(\Cn)\}$. 

One has $\alpha_n^*(\omega_i)=\zeta_n^i\omega_i$ and $\iota^*(\omega_i)=-\omega_i$.

If $n$ is odd, the local action of $\alpha_n$ near the fixed points $P_1$, $P_{-1}$, $P_{\infty}$ is $\zeta_n$, $\zeta_n$ and $\zeta_n^{\frac{n-1}{2}}$ respectively. 

If $n$ is even, the local action of $\alpha_n$ near the fixed points $P_1$, $P_{-1}$, is $\zeta_n$. The square $\alpha_n^2$ fixes also the points $P_{\infty}^{1}$ and $P_{\infty}^{2}$ and its the local action near these points is $\zeta_n^2$. 
\end{proposition}

\proof Since $\Cn$ is a hyperelliptic curve, the genus of $\Cn$ and the basis of $H^{1,0}(\Cn)$ can be computed as in \cite[Chapter 1 Section 2]{ACGH}. The explicit equations of $\alpha_n$ and $\iota_n$ allows one to check the number of the points with non trivial stabilizer. The local action of $\alpha_n$ near the points with non trivial stabilizer follows by \cite[Theorem 7]{H}, see \cite[Example 2.5]{GP} for a more the direct computation in the case of the curve $\Cn$. One needs to observe that if $n$ is odd then $n\frac{n-1}{2}\equiv 1\mod n$.\endproof

\subsection{The varieties $\Ynl$}\label{subsect: The varieties Ynl} Let $l$ be an integer bigger than 1. Let us denote by $\left(\Cn\right)^{\times l}$ the product of $\Cn$ by itself $l$ times . 
\begin{definition}
	
Let $\Gnl$ be the subgroup of $\Aut\left(\left(\Cn\right)^{\times l}\right)$ defined by  
$$\Gnl:=\left\{\alpha_n^{a_1}\times \alpha_{n}^{a_2}\times\ldots\times\alpha_{n}^{a_l}\mbox{ such that }\sum_{i=1}^la_i\equiv 0\mod n\right\}$$
and $\Ynl$ be the variety
 $\Ynl:=\left(\Cn\right)^{\times l}/\Gnl.$

We define $Y_{(n)}^{(1)}=\Cn$ and denote by $\pi_{(n)}^{(l)}$ the quotient map $\pi_{(n)}^{(l)}:\left(\Cn\right)^{\times l}\ra \Ynl.$
\end{definition}

We observe that $\Gnl\simeq(\Z/n\Z)^{l-1}$ and $$\Gnl=\langle\alpha_n\times\alpha_n^{n-1}\times \id\times\ldots\id, \alpha_n\times\id\times \alpha_n^{n-1}\times \id\times\ldots\times\id, \ldots, \alpha_n\times\id\times\ldots\times\id\times\alpha_n^{n-1}\rangle.$$
\begin{proposition}
	The Hodge numbers $h^{i,0}$ of any desingularization of the singular variety $\Ynl$ are $h^{0,0}=1$, $h^{l,0}=g(\Cn)$ and $h^{i,0}=0$ for any $i\neq0,l$.
\end{proposition}
\proof The Hodge numbers  $h^{i,0}$ are birational invariants, and they are  $$h^{i,0}(\Ynl)=\dim\left(H^{i,0}(\left(\Cn\right)^{\times l})\right)^{\Gnl}.$$ The result follows by Proposition \ref{prop: properties curve Cn}. \endproof

\begin{definition}
The automorphism $\alpha_{\Ynl}\in\Aut(\Ynl)$ is the automorphism induced by $\alpha_n\times \id\times\ldots\times\id\in\Aut\left(\left(\Cn\right)^{\times l}\right)$ on $\Ynl$.

The automorphism $\iota_{\Ynl}\in\Aut(\Ynl)$ is the automorphism induced by $\iota_n\times \id\times\ldots\times\id\in\Aut(\left(\Cn\right)^{\times l})$ on $\Ynl$.
\end{definition}

Let $K^{(l)}_{(n)}$ be the group $\left\{\alpha_n^{a_1}\times \ldots\times\alpha_n^{a_l},\ a_i\in\Z\right\}$. Clearly $K_{(n)}^{(l)}/G_{(n)}^{(l)}\simeq \Z/n\Z$ and a generator is induced by $\alpha_n\times\id \ldots\times\id\in K_{(n)}^{(l)}$. So $\Ynl/\alpha_{\Ynl}$ is birational to $\Cn^{\times l}/K_{(n)}^{(l)}\simeq \left(\mathbb{P}^1\right)^{\times l}$. More precisely,
denoted by $u_j$, $v_j$ and $w_j$ the variables of the $j$-th copy of $\Cn$ in $\left(C_{(n)}\right)^{\times l}$, following proposition holds.
\begin{proposition}
	The variety $\Ynl$ is a (singular) $n:1$ cover of $\left(\mathbb{P}^1\right)^{\times l}$ whose cover automorphism is $\alpha_{\Ynl}$ and whose equation is 
	\begin{equation}\label{eq: n cover of P1P1P1...}U^n=\prod_{i=1}^l(v_i-w_i)(v_i+w_i)w_i^{n-2}.\end{equation}
\end{proposition} 
\proof 
The functions $U:=\prod_{i=1}^lu_i$, $v_i$, $w_i$, $i=1,\ldots,l$ of $\Cn^{\times l}$ are invariant for $G_{(n)}^{(l)}$ and they satisfy \eqref{eq: n cover of P1P1P1...}. So $\Ynl$ is an $r:1$ cover of the singular variety defined by \eqref{eq: n cover of P1P1P1...}, where $r$ is a positive integer, possibly 1. Since the $n^l:1$ cover $\Cn^{\times l}\ra \left(\Cn/\alpha_n\right)^{\times l}=\left(\mathbb{P}^1\right)^{\times l}$ factorizes through the $n^{l-1}:1$ cover $\pi_{(n)}^{(l)}:\Cn^{\times l}\ra \Ynl$, one obtains that $\Ynl$ is an $n:1$ cover of $\left(\mathbb{P}^1\right)^{\times l}$. Therefore, \eqref{eq: n cover of P1P1P1...} is an equation of $\Ynl$, i.e. $r=1$. Since the automorphism $\alpha_{\Ynl}$ acts trivially on $v_i$ and $w_i$, and $\alpha_{\Ynl}(U)=\zeta_n U$, $\alpha_{\Ynl}$ is the cover automorphism. \endproof

\subsection{Relations between the varieties $\Ynl$ and $Y_{(m)}^{(h)}$}\label{subsec: relations between Ynl Ymh}
In this section we show that there are strong interactions between the $\Ynl$ and $Y_{(m)}^{(h)}$ if there are relations between the numbers $n,m,l,h$. In particular we show that:\begin{itemize}
	\item   $\Ynl$ is a $k:1$ cover of $Y_{(n/k)}^{(l)}$ (Proposition \ref{prop: covers between Ynl and Y(n/k)l});
	\item  $\Ynl$ is an $n:1$ quotient of $Y_{(n)}^{(l-h)}\times Y_{(n)}^{(h)}$ (Proposition \ref{prop: Ynl quotients of products of Ynl})
	\item there is an isotrivial fibration $\Ynl\ra\left(\mathbb{P}^1\right)^{\times (l-h)}$ whose fibers are isomorphic to $Y_{(n)}^{(h)}$ (Proposition \ref{prop: isotrivial firbation on Ynl}). 
	\end{itemize}

\begin{proposition}\label{prop: covers between Ynl and Y(n/k)l}
	Let $m$ and $n$ be two positive integers such that $m|n$ and $k$ be the integer such that $n=km$.
	
	Then $\Ynl$ is a $k:1$ cover of $Y_{(m)}^{(l)}$ and the cover automorphism is $\alpha_{\Ynl}^{m}$.
\end{proposition}
\proof If $l=1$ we are just considering the curves $C_{(n)}$ and $C_{(m)}$. Given $\Cn$ and $\alpha_n$ as in \eqref{eq: equation Cn} and \eqref{eq: equation alphan}, one obtains that $u':=u^k$, $v':=v$, $w':=w$ are invariant functions of the order $k$ automorphism $\alpha_n^m$ which satisfy the equation $u'^m=(v'-w')(v'+w')w'^{n-2}$. The automorphism $\alpha_n$ acts on $u'=u^k$ as the multiplication by $\zeta_m$. So the quotient $\Cn/\alpha_n^m$ is isomorphic to $C_{(m)}$ and the automorphism $\alpha_n$ induces on $C_{(m)}=\Cn/\alpha_n^m$ the automorphism $\alpha_m$. We observe that $\alpha_m$ is induced equivalently by $\alpha_n$ and $\alpha_n^k$ (since we are considering the quotient by $\alpha_n^m$). 
If $l>1$, then the cover is induced by the following commutative diagram:
$$\xymatrix{\left(C_{(n)}\right)^{\times l}\ar[rr]^-{n^{l-1}:1}\ar[d]_{k^l=(n/m)^l:1}&& Y_{(n)}^{(l)}:=\left(\left(\Cn\right)^{\times l}\right)/ G_{(n)}^{(l)}\ar[d]^{k:1}\\
\left(C_{(m)}\right)^{\times l}\ar[rr]^-{m^{l-1}:1}&& Y_{(m)}^{(l)}:=\left(\left(C_{(m)}\right)^{\times l}\right)/ G_{(m)}^{(l)}
}$$
More explicitly, the equation of $\Ynl$ and of $Y_{(m)}^{(l)}$ are given in \eqref{eq: n cover of P1P1P1...}: $\Ynl$ is
$U^n=\prod_{i=1}^l(v_i-w_i)(v_i+w_i)w_i^{n-2}$ where $U:=\prod_{i=1}^l u_i$  and  $Y_{(m)}^{(l)}$ is  $U'^m=\prod_{i=1}^l(v_i'-w_i')(v_i'+w_i')w_i'^{n-2}$  where $U':=\prod_{i=1}^l u_i'$ and $u_i'=u_i^k$. So $U'=U^k$ and the map $\Ynl\ra Y_{(m)}^{(l)}$ is the quotient map by $\alpha_{Y_{(n)}^{(l)}}^{m}$.\endproof

Let $l$ and $h$ be two positive integers such that $h<l$.
The group $G_{(n)}^{(l-h)}$ can be embedded in $G_{(n)}^{(l)}$ as the normal subgroup $$G_{(n)}^{(l-h)}:=\left\{\alpha_n^{a_1}\times \alpha_{n}^{a_2}\times\ldots\times\alpha_{n}^{a_l}\mbox{ such that }\sum_{i=1}^{l-h}a_i\equiv 0\mod n,\ \ a_j=0\mbox{ for any }j>l-h\right\}.$$ We also define $$H_{(n)}^{(h)}:=\left\{\alpha_n^{a_1}\times \alpha_{n}^{a_2}\times\ldots\times\alpha_{n}^{a_l}\mbox{ such that }\sum_{i=l-h+1}^{l}a_i\equiv 0\mod n,\ \ a_j=0\mbox{ for any }j<l-h+1\right\}.$$

	\begin{proposition}\label{prop: Ynl quotients of products of Ynl}
		Let $l$ and $h$ be two positive integers such that $l>h$. Then $\Ynl$ is the quotient of $Y_{(n)}^{(l-h)}\times Y_{(n)}^{(h)}$ by the order $n$ automorphism $\alpha_{Y_{(n)}^{(l-h)}}\times\alpha_{Y_{(n)}^{(h)}}^{n-1}$.
	\end{proposition}
\proof Let $G_{(n)}^{(l-h)}$ and $H^{(h)}_{(n)}$ be as above. We observe that $H_{(n)}^{(h)}\simeq (\Z/n\Z)^{h-1}\simeq G_{(n)}^{(h)}$ and that $G_{(n)}^{(l-h)}\times H_{(n)}^{(h)}\simeq (\Z/n\Z)^{l-2}$ is a normal subgroup of $\Gnl$. Moreover, $\Gnl/\left(G_{(n)}^{(l-h)}\times H_{(n)}^{(h)}\right)$ is cyclic and generated by the order $n$ automorphism $\overline{\alpha}:=\alpha_n^{a_1}\times \alpha_{n}^{a_2}\times\ldots\times\alpha_{n}^{a_l}$ with $a_1=1$, $a_{l-h+1}=n-1$ and $a_j=0$ for any other $j$.

Since $$\begin{array}{l}\left(\Cn\right)^{\times l}/\left(G_{(n)}^{(l-h)}\times H_{(n)}^{(h)}\right)=\left(\left(\Cn\right)^{\times (l-h)}\times \left(\Cn\right)^h\right)/\left(G_{(n)}^{(l-h)}\times H_{(n)}^{(h)}\right)=\\
=\left(\left(\left(\Cn\right)^{\times (l-h)}\right)/G_{(n)}^{(l-h)}\right)\times\left( \left(\left(\Cn \right)^{\times h}\right)/H_{(n)}^{(h)}\right)=Y_{(n)}^{(l-h)}\times Y_{(n)}^{(h)}\end{array}$$ one obtains $$\Ynl\simeq \left(Y_{(n)}^{(l-h)}\times Y_{(n)}^{(h)}\right)/\left(\alpha_{Y_{(n)}^{(l-h)}}\times\alpha_{Y_{(n)}^{(h)}}^{n-1}\right),$$ where  $\left(\alpha_{Y_{(n)}^{(l-1)}}\times\alpha_{Y_{(n)}^{(h)}}^{n-1}\right)$ is induced on $Y_{(n)}^{l-h}\times Y_{(n)}^{h}$ by $\overline{\alpha}$.\endproof

Let us consider the case $h=l-1$ in Proposition \ref{prop: Ynl quotients of products of Ynl}. We obtain the following diagram:
$$\xymatrix{\left(\Cn \right)^{\times l}\ar[rrr]^{/\Gnl}\ar[d]_{/G_{(n)}^{(l-1)}}&&&\Ynl\\
Y_{(n)}^{(l-1)}\times C_n\ar[urrr]_{\alpha_{Y_{(n)}^{(l-1)}}\times\alpha_{n}^{n-1}}}$$
which shows that $\Ynl$ is the cyclic quotient (of order $n$) of $Y_{(n)}^{(l-1)}\times\Cn$ by the automorphism $\alpha_{Y_{(n)}^{(l-1)}}\times\alpha_n^{n-1}$.
This provide an inductive construction of the varieties $Y_{(n)}^{(l)}$ from the varieties $Y_{(n)}^{(l-1)}$, which will be used in the sections \ref{sec: 3-fold} and \ref{sec: l-fold} to find an equation of $\Ynl$.

The construction of $\Ynl$ gives also fibrations (surjective maps with connected fibers) from $\Ynl$ to $\left(\mathbb{P}^1\right)^{\times (l-h)}$ which are induced by the projections $\left( \Cn\right)^{\times l}\ra \left( \Cn\right)^{\times (l-h)}$
and are described in the following proposition.

\begin{proposition}\label{prop: isotrivial firbation on Ynl}
Let $l$ and $h$ be two positive integers such that $l>h$. Then there is an isotrivial fibration $\Ynl\ra\left(\mathbb{P}^1\right)^{\times (l-h)}$ whose fibers are isomorphic to $Y_{(n)}^{(h)}$.
	\end{proposition}
\proof 
 We consider the action of $G_{(n)}^{(l)}\simeq (\Z/n\Z)^{l-1}$ on $\left(\Cn\right)^{\times l}$, the action of $\alpha_n$ on $\Cn$ and the projection on the first $(l-h)$ factors $p:\Cn^{\times l}\ra \Cn^{\times (l-h)}$.
Then we have the following diagram
$$\xymatrix{\left(\Cn\right)^{\times l}\ar[r]^p\ar[d]_{n^{l-1}:1}& \left(\Cn \right)^{\times (l-h)}\ar[d]^{n^{l-h}:1}\\
\Ynl:=\left( \Cn \right)^{\times l}/\Gnl\ar[r]^-{\mathcal{F}}&\left(\Cn/\alpha_n\right)^{\times (l-h)}\simeq \left(\mathbb{P}^1\right)^{\times (l-h)}
}$$
One can choose a set of generators of $\Gnl$ to be of the form $\alpha\times id \times\ldots\times\id\times \alpha_n^{l-1}$, $\id\times  \alpha\times \ldots\times\id\times \alpha_n^{l-1}$, $\ldots$, $\id\times\id\times\ldots\times \alpha\times\alpha^{l-1}$. So, the generic fiber of $\mathcal{F}$ is isomorphic to $\left(\Cn \right)^{\times h}/H_{(n)}^{(h)}$ and hence to $Y_{(n)}^{(h)}$.\endproof

In particular, if $h=l-1$, the fibration $\mathcal{F}$ is an isotrivial fibration with smooth fibers isomorphic to $\Cn$.

Given a multi-index $I=(i_1,\ldots, i_h)$ of length $h$ such that $1\leq i_1<i_2<\ldots<i_h\leq l$, we denote by $\mathcal{F}_{I}$ the fibration $\Ynl\ra\left(\mathbb{P}^1\right)^{l-h}$ induced by the projection of $\Cn\times\ldots\times \Cn$ on the factors in position $(i_1,\ldots, i_h)$.

In the particular case $l=2$, this gives two isotrivial fibrations $\mathcal{F}_i:\Yns\ra\mathbb{P}^1$ for $i=1,2$ whose generic fibers are isomorphic to $\Cn$.
\subsection{Kodaira dimensions}
\begin{proposition}\label{prop: Kodaira dimension}
It holds $k\left(Y_{(n)}^{(l)}\right)\leq k\left(Y_{(n)}^{(l-1)}\right)+1$.	
	\end{proposition}
\proof Let $\Xns$ be the minimal resolution of $\Yns$. The singularities of $\Yns$ are all contained in 3 singular fibers of the isotrivial fibration of $\Yns\ra \mathbb{P}^1\simeq \Cn/\alpha_n$ described in Proposition \ref{prop: isotrivial firbation on Ynl}. So $\Xns$ is endowed with an isotrivial fibration $\Xns\ra\mathbb{P}^1$, which coincides with $\Yns\ra\mathbb{P}^1$ outside the fibers over $1$, $-1$, $\infty$. The automorphism $\alpha_Y$ lifts to an automorphism $\alpha_X$ of the smooth model $\Xns$. The singular quotient variety  $\left(\Xns\times \Cn\right)/\left(\alpha_X\times\alpha_n^{n-1}\right)$ is birational to  $\left(\Yns\times \Cn\right)/\left(\alpha_Y\times\alpha_n^{n-1}\right)$, and thus, by Proposition \ref{prop: Ynl quotients of products of Ynl}, to $\Ynt$. The singularities of $\left(\Xns\times \Cn\right)/\left(\alpha_X\times\alpha_n^{n-1}\right)$ are all contained in the three fibers over $1$, $-1$ and $\infty$ of the isotrivial fibration $\left(\Xns\times \Cn\right)/\left(\alpha_X\times\alpha_n^{n-1}\right)\ra\mathbb{P}^1\simeq \Cn/\alpha_n^{n-1}$, whose generic fibers are isomorphic to $\Xns$. To resolve these singularities one introduces some divisors contained in these three singular fibers and one obtaines a variety, $\Xnt$, which is smooth and admits an isotrivial fibration on $\mathbb{P}^1$ whose generic fiber is isomorphic to $\Xns$. The automorphism $\alpha_X$ on $\Xns$ induces an automorphism of the quotient $\left(\Xns\times \Cn\right)/\left(\alpha_X\times\alpha_n^{n-1}\right)$ which lifts to an automorphism, still denoted by $\alpha_X$, on its desingularization $\Xnt$. By construction $\Xnt$ is birational to $\Ynt$ and the automorphism $\alpha_X$ is the one induced by $\alpha_Y$. One can iterate the previous process to obtain a smooth variety $\Xnl$, birational to $\Ynl$, which is a resolution of $\left(X_{(n)}^{(l-1)}\times \Cn\right)/\left(\alpha_X\times\alpha_n^{n-1}\right)$ such that the automorphism induced by $\alpha_X\times \id$ on $\left(X_{(n)^{(l-1)}}\times \Cn\right)/\left(\alpha_X\times\alpha_n^{n-1}\right)$ lifts to an automorphism $\Xnl$ denoted by $\alpha_X\in \Aut(\Xnl)$. Moreover, $\Xnl$ is endowed with an isotrivial fibration $\Xnl\ra\mathbb{P}^1$ whose generic fibers are isomorphic to $X_{(n)}^{(l-1)}$ (which is a smooth variety birational to $Y_{(n)}^{(l-1)}$). 	
By applying the ``easy addiction formula" (see a.g. \cite[Theorem 6.122]{U}) to the isotrivial fibration $\Xnl\ra\mathbb{P}^1$ with generic fiber isomorphic to $X_{(n)}^{(l-1)}$, one has $$k(\Xnl)\leq k\left(X_{(n)}^{(l-1)}\right)+\dim(\mathbb{P}^1).$$
Since the Kodaira dimension is a birational invariant, we have the following chain of (in)equalities which proves the statement $$k(\Ynl)=k\left(\Xnl\right)\leq k\left(X_{n}^{(l-1)}\right)+1=k\left(Y_{(n)}^{(l-1)}\right)+1.$$  \endproof
\begin{corollary} If $\Ynl$ is not of general type for a certain $l$, then $Y_{(n)}^{(h)}$ is not of general type for every $h>l$.
\end{corollary}
\proof By Proposition \ref{prop: Kodaira dimension}, $k\left(Y_{(n)}^{(h)}\right)\leq k\left(\Ynl\right)+(h-l)$. So if $\Ynl$ is not of general type, i.e. $k\left(\Ynl\right)<l$, then $k\left(Y_{(n)}^{(h)}\right)<l+(h-l)=h$.\endproof
We will give a more restrcitive bound on the Kodaira dimension of $\Ynl$ in Section \ref{sec: l-fold}.

\section{The surfaces $\Yns$ and $\Xns$}\label{sec: Yns}
In this section we consider the singular surface $\Yns$ and its minimal model. 
If $n>2$, we will prove that the smooth minimal model of $\Yns$ coincides with the minimal resolution $\Xns$ of $\Yns$. We prove that $\Xns$ admits an elliptic fibration, obtained by a base change (of order $n$) on an elliptic fibration on a rational surface, and that its Kodaira dimension is 1 if $n>4$ (if otherwise $n=3,4$, $\Xns$ is a K3 surface and so its Kodaira dimension is 0), and we compute its Hodge diamond and its Picard number (proving that it is an extremal elliptic surface).
If $n=2$, the surface $\Yns$ is rational, and we consider a birational model of $\Yns$ which is not a minimal surface but still admits a relatively minimal elliptic fibration, obtained by a base change (of order 2) from another rational  elliptic surface.

The results of this section are summarized by the following theorem.
 
\begin{theorem}\label{theorem surfaces} Let $n>2$ and $X_{(n)}^{(2)}$ be the minimal resolution of $Y_{(n)}^{(2)}$. Then:
	\begin{itemize}
		\item $\Xns$ is a minimal surface;
		\item $h^{1,0}\left(\Xns\right)=0$, $h^{2,0}\left(\Xns\right)=g\left(C_{(n)}\right)$ and if $n=3,4$, then $\Xns$ is a K3 surface, otherwise $k\left(\Xns\right)=1$;
		\item $$h^{1,1}\left(\Xns\right)=\rho\left(\Xns\right)=\left\{\begin{array}{ll}5n+5&\mbox{if  } n\equiv 1 \mod 2\\5n&\mbox{if  } n\equiv 0 \mod 2;\end{array}\right.$$  
		\item $\Xns$ is an elliptic modular surface.
	\end{itemize}
\end{theorem}

\subsection{The surfaces $\Yns$ as base change of a rational elliptic surface}\label{subsec: Yns as base change}
A rational elliptic surface can be constructed as a blow up of $\mathbb{P}^1\times\mathbb{P}^1$ in the base locus of a pencil $\mathcal{P}$ of generically smooth genus 1 curves, so of generically smooth curves of bidegree $(2,2)$.

We consider $\mathbb{P}^1\times \mathbb{P}^1$ with coordinates $((v_1:w_1),(v_2:w_2))$ and we put $$C_0:=V((v_1-w_1)(v_1+w_1)(v_2-w_2)(v_2+w_2)),\ C_{\infty}:=V(w_1^2w_2^2).$$ 
For a generic $t$, curve $C_{t}:=C_0+tC_{\infty}$ is smooth and so the pencil
\begin{equation}\label{eq: pencil P}\left(v_1-w_1\right)\left(v_1+w_1\right)\left(v_2-w_2\right)\left(v_2+w_2\right)+t\left(w_1^2w_2^2\right)\end{equation}
induces an elliptic fibration on the blow up of $\mathbb{P}^1\times\mathbb{P}^1$ in the points $C_0\cap C_{\infty}$, whose fiber over $\overline{t}$ is the strict transform $\widetilde{C_{\overline{t}}}$.
\begin{definition}\label{defi R}
We call $\mathcal{P}$ the pencil of $(2,2)$ curves generated by $C_0$ and $C_{\infty}$ with equation \eqref{eq: pencil P}, $R$ the surface obtained by blowing up the base points of $\mathcal{P}$ and $\mathcal{E}_R:R\ra\mathbb{P}^1_t$ the induced elliptic fibration.
\end{definition}

\begin{lemma}\label{lemma: R}
	The rational elliptic surface $R$ is a modular elliptic surface with Mordell--Weil group $MW(\mathcal{E}_R)=\Z/4\Z$ and singular fibers $I_1^*+I_4+I_1$. 
	\end{lemma}
\proof The singular members of the pencil $\mathcal{P}$ are $C_0$, $C_{\infty}$ and $C_1$. The latter is an irreducible curve singular in the point $((0:1),(0:1))$.
The base points of the pencil are the following, each with multiplicity 2: $$(p,q)\in\{((1:0),(1:1)),\ ((1:0),(-1:1)),\ ((1:1),(1:0)),\ ((-1:1),(1:0))\}.$$ Blowing up each of these points twice one introduces a $(-2)$- curve and a $(-1)$-curve for each of them. We will denote by $E_{(p,q)}^{(1)}$ and $E_{(p,q)}^{(2)}$ the two exceptional divisors over $(p,q)$, such that $\left(E_{(p,q)}^{(h)}\right)^2=-h$ and $E_{(p,q)}^{(1)}E_{(p,q)}^{(2)}=1$.
Then the canonical divisor of $R$ is $$K_R:=-2h_1-2h_2+\sum_{j=1}^2j\left(E_{((1:0),(1:1))}^{(j)}+E_{((1:0),(-1:1))}^{(j)}+E_{((1:1),(1:0))}^{(j)}+E_{((-1:1),(1:0))}^{(j)}\right)$$ The linear system of the divisor $F_R:=-K_R$ defines the elliptic fibration $\mathcal{E}_R:R\ra\mathbb{P}^1_t$. The exceptional curves $E_{((1:0),(1:1))}^{(1)}$, $E_{((1:0),(-1:1))}^{(1)}$, $E_{((1:1),(1:0))}^{(1)}$, $E_{((-1:1),(1:0))}^{(1)}$
are $(-1)$-curves and sections of the fibration. Denoted by $\widetilde{E_{p\times \mathbb{P}^1}}$ and 
 $\widetilde{E_{\mathbb{P}^1\times q}}$ the strict transform of the curves $p\times \mathbb{P}^1$ and $\mathbb{P}^1\times q$ respectively, the divisor
$$
\widetilde{E_{(1:1)\times \mathbb{P}^1}}+
\widetilde{E_{\mathbb{P}^1\times(1:1)}}+
\widetilde{E_{(-1:1)\times\mathbb{P}^1}}+
\widetilde{E_{\mathbb{P}^1\times(-1:1)}}$$
is orthogonal to $F_R$ (and indeed linearly equivalent to $F_R$) and is supported on configuration of curves which corresponds to a fiber of type $I_4$. Similarly the divisor
$$2(\widetilde{E_{\mathbb{P}^1\times(1:0)}}+\widetilde{E_{(1:0)\times \mathbb{P}^1}})+E_{((1:0),(1:1))}^{(2)}+E_{((1:0),(-1:1))}^{(2)}+E_{((1:1),(1:0))}^{(2)}+E_{((-1:1),(1:0))}^{(2)}$$
is supported on a configuration of curves which corresponds to a fiber of type $I_1^*$. Therefore the fibers of $\mathcal{E}_R$ over $t=0$ and $t=\infty$ are of type $I_4$ and $I_1^*$ respectively.  Hence $R$ is an extremal rational elliptic surface, with singular fibers $I_1^*+I_4+I_1$ and Mordell--Weil group $MW(\mathcal{E}_R)\simeq \Z/4\Z$.

The surface $R$ is not isotrivial (i.e. the $j$-function of $\mathcal{E}_R:R\ra\mathbb{P}^1$ is non constant since there are singular fibers of type $I_n$); does not admit reducible fibers of type $II^*$ and $III^*$; is extremal. Thus by \cite[Theorem 3.6]{N}, $R$ is a modular surface.\endproof

\begin{proposition}\label{prop: Yns birat to base change}
The surface $\Yns$ is birational to a base change of order $n$ of $R$, branched over the reducible fibers, and thus to a surface with an elliptic fibration.
\end{proposition}
\proof  Let $b:\mathbb{P}^1_\tau\ra\mathbb{P}^1_t$ be the $n:1$ cover of $\mathbb{P}^1$ totally branched over 0 and $\infty$. So $b$ is simply given by $\tau\mapsto t:=\tau^n$ and the cover automorphism is $\tau\mapsto\zeta_n \tau$. We consider the fiber product:
$$\xymatrix{R\times_g\mathbb{P}^1\ar[r]^{n:1}\ar[d]_{\mathcal{E}}&R\ar[d]^{\mathcal{E}_R}\\
	\mathbb{P}^1_{\tau}\ar[r]^{b}&\mathbb{P}^1_t}$$
The, possibly singular, surface $R\times_b\mathbb{P}^1$ naturally arise with the map $\mathcal{E}$ to $\mathbb{P}^1_{\tau}$ and there exists a birational model of $R\times_b\mathbb{P}^1$ which is a smooth surface with an elliptic fibration induced by $\mathcal{E}_R$.

By \eqref{eq: pencil P} a birational model of $R$ has the following equation:
$$t=-\frac{(v_1-w_1)(v_1+w_1)(v_2-w_2)(v_2+w_2)}{(w_1w_2)^2}.$$
So a birational model of $R\times_b\mathbb{P}^1$ has the following equation 
$$\tau^n=-\frac{(v_1-w_1)(v_1+w_1)(v_2-w_2)(v_2+w_2)}{(w_1w_2)^2}.$$
After a change of sign and the birational transformation $\tau\mapsto U/(w_1w_2)$, one recognizes the equation \eqref{eq: n cover of P1P1P1...}
with $l=2$. 
\endproof

Let us denote by $Z_{(n)}$ the smooth surface resolving the singularities of $R\times\mathbb{P}^1$, so that $Z_{(n)}$ admits the relatively minimal elliptic fibration $\mathcal{E}:Z_{(n)}\ra\mathbb{P}^1$ induced by $\mathcal{E}_R:R \ra\mathbb{P}^1$.

\begin{proposition}\label{prop: the surface Z_n}
Let $\mathcal{E}:Z_{(n)}\ra\mathbb{P}^1$ be the elliptic fibration induced by $\mathcal{E}_R:R\ra\mathbb{P}^1$ by the base change of order $n$ branched over the reducible fibers $I_1^*+I_4$. Then $$h^{1,0}(Z_{(n)})=0,\ h^{2,0}(Z_{(n)})=g(\Cn),\ h^{1,1}(Z_{(n)})=10(g(\Cn)+1).$$ The Euler characteristic of $Z_{(n)}$ is $ \ e(Z_{(n)})=12(g(\Cn)+1)$. 

If $n\neq 2$, $Z_{(n)}$ is a minimal surface and 
$$Z_{(n)}\mbox{ is }\left\{\begin{array}{lll}\mbox{ a rational surface }&\mbox{ if }n=2&\mbox{ so }k\left(Z_{(n)}\right)=\infty\\
\mbox{ a K3 surface }&\mbox{ if }n=3,4&\mbox{ so }k\left(Z_{(n)}\right)=0\\
\mbox{ a properly elliptic surface }&\mbox{ if }n>4&\mbox{ so }k\left(Z_{(n)}\right)=1.\end{array} \right.$$

The Mordell--Weil group of $\mathcal{E}$ is $MW(\mathcal{E})=\Z/4\Z$ and the singular fibers are
$$
\left\{\begin{array}{lll}I_n^*+I_{4n}+nI_1&\mbox{ if }n\equiv 1\mod 2\\
	I_{n}+I_{4n}+nI_1&\mbox{ if }n\equiv 0\mod 2.\end{array} \right.$$
\end{proposition}
\proof 
In \cite[Table VI.4.1]{M} the effect of a base change of order $n$ branched on a singular fiber is described: if the branch fiber is of type $I_m$, the base change produces a fiber of type $I_{nm}$; if the branch fiber is of type $I_{m}^*$ and $n$ is odd, the base change produces a fiber of type $I_{nm}^*$; if the branch fiber is of type $I_m^*$ and $n$ is even, the base change produces a fiber of type $I_{nm}$. For each fiber of $\mathcal{E}_R:R\ra\mathbb{P}^1$ which is not a branch fiber, the base change produces $n$ copies of that fiber on $\mathcal{E}:Z_{(n)}\ra\mathbb{P}^1$ and this proves the statement on the singular fibers of $\mathcal{E}:Z_{(n)}\ra\mathbb{P}^1$. Due to the presence of these reducible fibers, $Z_{(n)}$ is not a fibration of product type, hence $h^{1,0}$ equals the genus of the base of the fibration, which is 0. The Euler characteristic of an elliptic fibration is the sum of the Euler characteristic of its singular fibers, and using that $e\left(I_n^*\right)=n+6$, $e\left(I_{m}\right)=m$, one computes the Euler characteristic of $Z_{(n)}$. Since $Z_{(n)}$ admits an elliptic fibration, $c_1^2(Z_{(n)})=0$. By Noether formula the holomorphic characteristic $\chi(Z_{(n)})=e(Z_{(n)})/12$. Hence $h^{2,0}(Z_{(n)})=\left(e(Z_{(n)})/12\right)-1$.
Since $e(Z_{(n)})=\sum_i(-1)^i b_i(Z_{(n)})$, where $b_i(Z_{(n)})$ are the Betti numbers, one obtains
$$2+h^{2,0}(Z_{(n)})+h^{0,2}(Z_{(n)})+h^{1,1}(Z_{(n)})=e(Z_{n}),\ \ i.e.\ h^{1,1}(Z_{(n)})=5e(Z_{(n)})/6$$ 

For each reducible fiber $F$ we denote by $r(F)$ the number of its irreducible components, so $r(I_m)=m$ and $r(I_m^*)=m+5$. By the Shioda--Tate formula, the Picard number of $Z_{(n)}$ is $$\rho(Z_{(n)})=2+\sum_{F\mbox{ red. fiber }}\left(r(F)-1\right)+\rk(MW(\mathcal{E})).$$ So if $n$ is odd, one obtains $\rho(Z_{(n)})=2+(4n-1)+(n+5-1)+\rk(MW(\mathcal{E}))$ and if $n$ is even one obtains $\rho(Z_{(n)})=2+(4n-1)+(n-1)+\rk(MW(\mathcal{E}))$. In both the cases one finds $\rho(Z_{(n)})=h^{1,1}(Z_{(n)})+\rk(MW(\mathcal{E}))$. This forces the $\rk(MW(\mathcal{E}))=0$, so the Mordell--Weil group is a (possibly trivial) torsion group. The (torsion) sections of the fibration $\mathcal{E}_R:R\ra\mathbb{P}^1$ induce sections of the fibration $\mathcal{E}:Z_{(n)}\ra\mathbb{P}^1$ (of the same order). So $\Z/4\Z\subset MW(\mathcal{E})$. If $n$ is odd, the presence of a fiber of type $I_n^*$ forces the torsion part of the Mordell--Weil group to be contained in $\Z/4\Z$. If $n$ is even one needs to use the height formula to conclude the statement on the Mordell--Weil group. The definition and the details on this formula can be found e.g. in \cite{SS}. Let $P$ be the 4-torsion section of $\mathcal{E}$, and $2P:= P+_{MW}P$, where $+_{MW}$ is the sum with respect to the group law in the Mordell--Weil group. So $2P$ is a 2-torsion section of $\mathcal{E}$. By the height formula any 2-torsion section of $\mathcal{E}$ intersects the reducible fibers exactly as $2P$ does. This implies that a 2-torsion section of $\mathcal{E}$ is necessarily $2P$, and so $MW(\mathcal{E})$ is cyclic. Let $Q$ be a generator of $\MW(\mathcal{E})$. If $MW(\mathcal{E})\neq \Z/4\Z$, $Q$ is a torsion section whose order is strictly bigger than 4. By the height formula the sum of contributions of the reducible fibers to the height of $Q$ must add up to $n$. The contribution due to a fiber of type $I_{n}$ is less than or equal to $n/4$. Moreover, since $Q$ has order strictly bigger than $P$, either $\contr_{I_n}(Q)<\contr_{I_n}(P)=n/4$ or $\contr_{I_{4n}}(Q)<\contr_{I_{4n}}(P)=3n/4$. Then $\contr_{I_n}(Q)+\contr_{I_{2n}}(Q)<n$ and we  conclude $MW(\mathcal{E})=\langle P \rangle$.

The classification of the surface $Z_{(n)}$ follows directly by \cite[Lemma III.4.6]{M}.\endproof

\begin{corollary}\label{cor: Zn modulare}
The surface $Z_{(n)}$ is a modular elliptic surface.
\end{corollary}
\proof By Proposition \ref{prop: the surface Z_n}, the elliptic fibration $\mathcal{E}:Z_{(n)}\ra\mathbb{P}^1$ satisfies the following conditions:
\begin{itemize}
	\item $\mathcal{E}$ is an extremal elliptic fibration with a section;
	\item $\mathcal{E}$ has no fibers of type $II^*$ and $III^*$;
	\item the $j$-invariant of $\mathcal{E}$ is non constant.
\end{itemize}
One concludes (as in Lemma \ref{lemma: R}) by \cite[Theorem 3.6]{N}.\endproof
\begin{corollary}\label{cor: weierstrass Xns}
The elliptic fibration $\mathcal{E}:Z_{(n)}\ra\mathbb{P}^1_{(\tau:\sigma)}$ admits the following Weierstrass equation:
\begin{equation}\label{eq: weierstrass Xns} y^2=x(x^2+\sigma^{\varepsilon}(\sigma^n-2\tau^n)x+\tau^{2n}\sigma^{2\varepsilon}).\end{equation}
where $\varepsilon=0$ if $n\equiv 0\mod 2$ and $\varepsilon=1$ if $n\equiv 1\mod 2$.
\end{corollary}
\proof The rational elliptic surface $R$ is rigid and its Weierstrass equation is determined (up to multiplication by constants) by the requirement that it has a fiber of type $I_1^*$ at infinity and a fiber of type $I_4$ over zero. So a Weierstrass equation for $R$ is $y^2=x(x^2+(1-2t)x+t^2)$. Applying the base change $\tau^n=t$ and considering the homogenous equation one finds the equation \eqref{eq: weierstrass Xns}.\endproof

The automorphism $\alpha_{\Yns}$ induces an automorphism on $Z_{(n)}$ which will be denoted still by $\alpha_{\Yns}$, with a slight abuse of notation. Considering the non-homogeneous equation obtained by putting $\sigma=1$ in \eqref{eq: weierstrass Xns}, the action of $\alpha_{\Yns}$ is just \begin{equation}\label{eq: action alpha surface tau}(\tau,x,y)\mapsto(\zeta_n \tau,x,y).\end{equation}

\subsection{The geometry of the minimal resolution $\Xns$ of $\Yns$}\label{subsec: geometry of minimal resolution surface}
In this section we assume $n>2$ and we describe the minimal resolution $\Xns$ of $\Yns$. It is naturally endowed with an elliptic fibration $\mathcal{E}:\Xns\ra\mathbb{P}^1$ and we show that it is the minimal model of $\Yns$. In particular this implies that $\Xns $ coincides with $Z_{(n)}$ and allows us to give a basis of the N\'eron--Severi group supported on irreducible components of reducible fibers and on sections of the fibration. 
To construct the minimal resolution of $\Yns$ we first identify its singularities.

\begin{lemma}\label{sub: sing of Yns}\label{prop: sing of Yns}
	If $n$ is odd, the surface $\Yns$ is singular in 9 points: it has 5 singularities of type $\frac{1}{n}\left(1,n-1\right)$ and 4 of type $\frac{1}{n}\left(1,\frac{n+1}{2}\right)$.
	
	If $n>2$ is even, the surface $\Yns$ is singular in 10 points: it has 4 singularities of type $\frac{1}{n}\left(1,n-1\right)$, 4 of type $\frac{1}{n/2}\left(1,n/2-1\right)$ and 2 of type $\frac{1}{n/2}\left(1,1\right)$.
\end{lemma}
\proof
{\bf If $n$ is odd}, the 9 points $(P_{\pm 1}\times P_{\pm 1})$, $(P_{\infty}\times P_{\pm 1})$, $(P_{\pm 1}\times P_{\infty})$, $(P_{\infty}\times P_{\infty})\in \Cn\times \Cn$ are fixed by $\alpha_n\times\alpha_n^{n-1}$ and thus their images for $\pi_{(n)}^{(2)}$ are singular. The type of singularities is determined by the local action of $\alpha_n\times\alpha_n^{n-1}$ near to these points, which can be found by Proposition \ref{prop: properties curve Cn}. In the following diagram 
we summarize the information on the singularities

$$
\begin{array}{|c|c|c|}
\hline point&\mbox{local
	action of }\alpha_n\times\alpha_n^{n-1} &\mbox{singularities of }\Yns\\
\hline (P_{\pm 1}\times P_{\pm 1})&(\zeta_n,\zeta_n^{n-1})&\frac{1}{n}(1,n-1)\\
\hline
(P_{\infty}\times P_{\infty})&(\zeta_n^{(n-1)/2},\zeta_n^{(n+1)/2})&\frac{1}{n}(1,n-1)\\
\hline
(P_{\pm 1}\times P_{\infty})&(\zeta_n,\zeta_n^{(n+1)/2})&\frac{1}{n}(1,(n+1)/2)\\
\hline
(P_{\infty},P_{\pm 1})&(\zeta_n^{(n-1)/2},\zeta_n^{n-1})&\frac{1}{n}(1,(n+1)/2)\\
\hline
\end{array}$$
(we used that $(n-1)^2/2\equiv_n (n+1)/2$ since 2 is invertible and $(n-1)^2\equiv_n 1\equiv_n n+1$).

{\bf If $n$ is even}, then $\Yns$ has $10$ singularities which are the image for $\pi_{(n)}^{(2)}$ of the points $(P_{\pm 1}\times P_{\pm 1})$, $(P_{\infty}^j\times P_{\pm 1})$, $(P_{\pm 1}\times P_{\infty}^j)$, $(P_{\infty}^j\times P_{\infty}^h)$, $j=1,2$, $h=1,2$, of $\Cn\times\Cn$. The points $(P_{\pm 1}\times P_{\pm 1})$ are fixed by $\alpha_n\times\alpha_n^{n-1}$, so each of them gives one canonical singularity in $\Yns$. The points $(P_{\infty}^j\times P_{1})$, $j=1,2$ form a unique orbit for $\alpha_n\times\alpha_n^{n-1}$, and then they correspond to a unique singularity, which is a quotient singularity of order $n/2$. Similarly the unique image of the points  $(P_{\infty}^j\times P_{-1})$ (resp. $(P_{1}\times P_{\infty}^j)$,  $(P_{-1}\times P_{\infty}^j)$) is a quotient singularity of order $n/2$ . The $4$ points $(P_{\infty}^j\times P_{\infty}^h)$, $j=1,2$, $h=1,2$ form $2$ orbits, each containing 2 points. This gives other 2 quotient singularities of order $n/2$. These two singularities are canonical (because, by definition, the action of $\alpha_n\times\alpha_n^{n-1}$ is such that the action on the first coordinate is the opposite of the action on the second one). The other 4 singularities of order $n/2$ are of type $\frac{1}{n/2}(1,1)$.\endproof

To construct $\Xns$ one has to resolve the singularities of $\Yns$. The desingularization of the singularities on surfaces introduces configurations of rational curves, which are well known and often called Hirzebruch--Jung strings. We list here just the information needed in our context (see e.g.\ \cite[III.5]{BHPV})
\begin{itemize}
\item The singularities of type  $\frac{1}{k}(1,k-1)$ are the  singularities of type $A_{k-1}$. To resolve them one introduces a chain of $(k-1)$ rational curves with self intersection $-2$;
\item to resolve the singularities of type $\frac{1}{k}(1,\frac{k+1}{2})$ one introduces two
rational curves meeting in a point. They have self intersections $-2$ and $-\frac{k+1}{2}$ respectively;
\item to resolve the singularities of type $\frac{1}{k}(1,1)$ one introduces one curve with self intersection $-k$.
\end{itemize}

We denote by $\pi:\Cn\times \Cn\ra \Yns$ the quotient map and by $\beta:\Xns\rightarrow\Yns$ the minimal resolution of the singularities of $\Yns$.
 
\subsubsection{The case $n\equiv 1\mod 2$}\label{subsubsec: desing surfaces odd n}
\begin{center}
	\begin{figure}[h]
		\includegraphics[width = 0.5\textwidth]{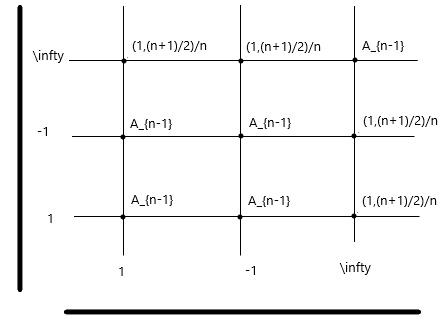}
		\caption{Singularities of $Y_{(n)}^{(2)}$.}
		\label{figure: sing}
	\end{figure}
\end{center}

Let $P_\epsilon\times P_\eta$ be one of the fixed points described in Lemma \ref{sub: sing of Yns} (see also Figure \ref{figure: sing}). Let us denote by $E_{P_\epsilon\times P_\eta}^{(j)}\subset \Xns$ the exceptional divisors such that $\beta(E_{P_\epsilon\times P_\eta}^{(j)})=\pi(P_\epsilon\times P_\eta)$, i.e.\ $E_{P_\epsilon\times P_\eta}^{(j)}$ are the exceptional divisors resolving the singularities on the singular points $\pi(P_\epsilon\times P_\eta)$. So, if $(\epsilon,\eta)\in\{ (\pm 1,\pm 1), (\infty,\infty)\}$, then $j=1,\ldots, n-1$ and $E_{P_\epsilon\times P_\eta}^{(j)}$, $j=1,\ldots, n-1$ form a configuration $A_{n-1}$ of $(-2)$-curves. If $(\epsilon,\eta)\in\{ (\pm 1,\infty), (\infty,\pm 1)\}$, then $j=1,2$, $\left(E_{P_\epsilon\times P_\eta}^{(1)}\right)^2=-2$, $\left(E_{P_\epsilon\times P_\eta}^{(2)}\right)^2=-\frac{n+1}{2}$ and $\left(E_{P_\epsilon\times P_\eta}^{(1)}\right)\left(E_{P_\epsilon\times P_\eta}^{(2)}\right)=1$. Moreover, we denote by 
$\widetilde{E}_{(P_{\epsilon}\times\Cn)}$ the strict transform (in
$X_{(n)}^{(2)}$) of $\pi(P_{\epsilon}\times \Cn)\subset \Yns$,
and, similarly, we denote by
$\widetilde{E}_{(\Cn\times P_{\eta})}$ the strict transform of $\pi(\Cn\times P_{\eta})$. The curves $\widetilde{E}_{(P_{\epsilon}\times \Cn)}$ with $\epsilon=\pm 1,\infty$ are components of the reducible fibers of the fibration $\mathcal{F}_1$ (see Proposition \ref{prop: isotrivial firbation on Ynl}); the curves $\widetilde{E}_{(\Cn\times P_{\eta})}$ with $\eta=\pm 1,\infty$ are components of the reducible fibers of the fibration $\mathcal{F}_2$. These components are studied in \cite{P}: they are called ``central component" of the reducible fibers of the isotrivial fibrations $\mathcal{F}_1$ and $\mathcal{F}_2$ respectively (cf.\ \cite[Theorem 2.3]{P}). By \cite[Proposition 2.8]{P}, one obtains that each of these 6 curves is a rational $-2$ curve.

\begin{center}
	\begin{figure}[h]
		\includegraphics[width = 0.5\textwidth]{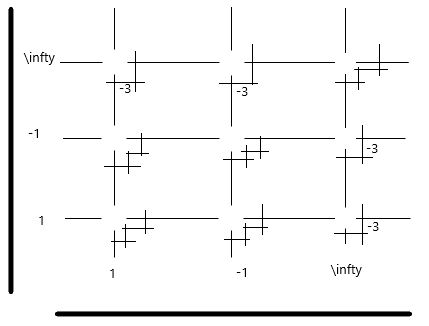}
		\caption{Minimal resolution of $Y_{(5)}^{(2)}$ (all the curves are $-2$ curves, with the exception of the 4 $(-3)$-curves).}
		\label{figure: X52}
	\end{figure}
\end{center}

\begin{center}
	\begin{figure}[h]
		\includegraphics[width = 0.5\textwidth]{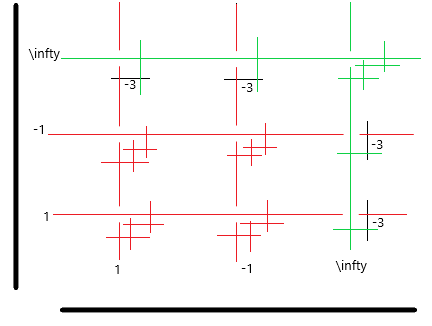}
		\caption{Minimal resolution of $Y_{(5)}^{(2)}$: the support of the divisor $F_1$ is red, the one of $F_2$ is green, the sections are black.}
		\label{figure: elliptic divisor X52}
	\end{figure}
\end{center}

This allows us to describe several rational curves in $\Xns$ (see Figure \ref{figure: X52} where the interesting curves in $\Xns$ are shown in case $n=5$). In particular we introduce two
interesting nef divisors:
\begin{align}\label{eq: F1F2} F_1:=\widetilde{E}_{(P_1\times \Cn)}+
\widetilde{E}_{(\Cn\times P_{1})}+
\widetilde{E}_{(P_{-1}\times \Cn)}+
\widetilde{E}_{(\Cn\times P_{-1})}+
\sum_{j=1}^{n-1} \left(E_{(1;1)}^{(j)}+E_{(-1;1)}^{(j)}+E_{(-1;-1)}^{(j)}+E_{(1;-1)}^{(j)}\right) \end{align}
$$F_2:=2\left( \widetilde{E}_{(P_{\infty}\times \Cn)}+ \widetilde{E}_{(\Cn\times P_{\infty})}+\sum_{j=1}^{n-1} E_{(\infty;\infty)}^{(j)} \right)+E_{(\infty;-1)}^{(1)}+E_{(\infty;1)}^{(1)}+E_{(1;\infty)}^{(1)}+E_{(-1;\infty)}^{(1)}.$$
The divisor $F_1$ is supported on $(-2)$-curves which form an $I_{4n}$ fiber and the divisor $F_2$ is supported on $(-2)$-curves which form an $I_{n}^*$ fiber (see Figure \ref{figure: elliptic divisor X52} for the case $n=5$). Moreover $F_1F_2=0$, so these are classes of fibers of the same fibration, which is necessarily a genus 1 fibration induced by the linear system $|F_1|$. In particular there are no $(-1)$-curves on $\Xns$ which are vertical with respect to $\varphi_{|F_1|}$. Moreover, the curves $E_{(\infty;-1)}^{(2)}$, $E_{(\infty;1)}^{(2)}$, $E_{(1;\infty)}^{(2)}$ and $E_{(-1;\infty)}^{(2)}$ are rational curves, they are sections of the fibration induced by $|F_1|$, and their self intersection is $\frac{n-1}{2}$. We recall that for any relatively minimal elliptic fibration with a rational basis, the self intersection of any section is equal to $-\left(h^{2,0}+1\right)$. Hence the self intersection of $E_{(\infty;-1)}^{(2)}$, $E_{(\infty;1)}^{(2)}$, $E_{(1;\infty)}^{(2)}$, $E_{(-1;\infty)}^{(2)}$ equals the self intersection of a section of $\mathcal{E}:Z_{(n)}\ra\mathbb{P}^1$. Thus, if there were $(-1)$-curves on $\Xns$, which have to be contracted to obtain the minimal model, they do not intersect these sections (otherwise the self intersection of the section would change after the contraction to obtain a minimal model). But the unique points which were blown up by $\beta$ are the singular points of $\Yns$, so there cannot be $(-1)$ curves which do not intersect at least one among the following curves: the fiber $F_1$; the fiber $F_2$; the sections $E_{(\infty;-1)}^{(2)}$, $E_{(\infty;1)}^{(2)}$, $E_{(1;\infty)}^{(2)}$ and $E_{(-1;\infty)}^{(2)}$. So, by comparison with Proposition \ref{prop: Yns birat to base change}, one obtains that $Z_{(n)}$ is isomorphic to $\Xns$ and the elliptic fibration $\mathcal{E}$ is the fibration induced by the linear system $|F_1|$.

\subsubsection{The case $n\equiv 0\mod 2$}\label{subsubsec: desing surfaces even n}
If {\bf $n$ is even}, one obtains the analogous result in a similar way, it suffices to slightly modify the divisor $F_2$ to obtain a fiber of type $I_{n}$. Indeed there are two singularities of type $A_{(n/2)-1}$ on $\Yns$. Let us denote by $E_{(\infty;\infty)^{(k)}}^{(j)}$, $k=1,2$, $j=1,\ldots, \left((n/2)-1\right)$ the $(n/2-1)$ irreducible curves resolving the $k$-th singularity of type $A_{(n/2)-1}$. They are $(-2)$-curves and the divisor 
\begin{align}\label{eq: F2 even}F_2:=\widetilde{E}_{(P_{\infty}\times\Cn)}+\sum_{j=1}^{n/2-1} E_{(\infty;\infty)^{(1)}}^{(j)}+\sum_{j=1}^{n/2-1} E_{(\infty;\infty)^{(2)}}^{(j)}+\widetilde{E}_{(\Cn\times P_{\infty})}
\end{align} is supported on a $I_{n}$-configuration of rational curves.
On $\Yns$ there are 4 singularities of type $\frac{1}{(n/2)}(1,1)$, which can be resolved by introducing one rational curve of self intersection $-n/2$. The four exceptional curves obtained by this resolution are 4 sections of the fibration (indeed the 4-torsion Mordell--Weil group).

\subsection{The proof of Theorem \ref{theorem surfaces}}

By Proposition \ref{prop: Yns birat to base change}, $\Yns$ is birational to a surface which is a base change of $R$ and we denoted by $Z_{(n)}$ the minimal elliptic surface resolving the singularities of the base change $R\times_{b}\mathbb{P}^1$.
In Sections \ref{subsubsec: desing surfaces odd n} and \ref{subsubsec: desing surfaces even n} it is proved that if $n>2$ the minimal resolution $\Xns$ of $\Yns$ is a minimal surface isomorphic to $Z_{(n)}$. In Proposition \ref{prop: the surface Z_n}, the Hodge diamond and the Kodaira dimension of $Z_{(n)}$ (and thus of $\Xns$) are computed and in Corollary \ref{cor: Zn modulare} it is proved that $Z_{(n)}$ (and thus $\Xns$) is a modular surface. This conclude the proof of the Theorem \ref{theorem surfaces}.

\begin{remark}{\rm In Proposition \ref{prop: covers between Ynl and Y(n/k)l} we observed that $\Ynl$ is a $k:1$ cover of $Y_{(n/k)}^{(l)}$ and the cover automorphism is $\alpha_{\Yns}^{k/n}$. Thanks to Theorem \ref{theorem surfaces} we can reinterpret this result  if $l=2$ in a more geometric way: the surface $\Xns$ with the elliptic fibration given in \eqref{eq: weierstrass Xns} is obtained by the elliptic fibration on $X_{(n/k)}^{(2)}$ by a base change of order $k$ branched on the fibers at 0 and $\infty$. The cover automorphism is $\alpha_{\Xns}^{k/n}$, where $\alpha_{\Xns}$ is induced by $\alpha_{\Yns}$ and its local action is given in \eqref{eq: action alpha surface tau}	}.
\end{remark}

\subsection{Translation by the torsion sections}
There are automorphisms of $\Cn^{\times l}$, not contained in the group $\Gnl$, which are indeed constructed by considering the automorphism $\iota_n$ on $\Cn$ and the permutations of the factors of the product $\Cn^{\times l}$. They induce automorphisms on $\Ynl$. In \cite{S} the quotients by some of these automorphisms are considered. Here, we prove that the actions induced by some of them on $\Xns$ preserve the elliptic fibration and in fact are translations by torsion sections. \\

We recall that $\iota_n$ is the hyperelliptic involution on $\Cn$ and denote by $\varepsilon$ the involution of $\Cn\times \Cn$ switching the two copies of $\Cn$. The automorphisms $\iota_n\times \iota_n$ and $\varepsilon\circ(\iota_n\times \id)$ of $\Cn\times \Cn$ induce automorphisms on $\Yns$
and on the smooth surface $\Xns$. In particular, we denote by $\iota_{\Xns}$ the involution of $\Xns$ induced by $\iota_n\times\iota_n$ and by $\gamma$ the order 4 automorphism induced on $\Xns$ by $\varepsilon\circ(\iota_n\times \id)$.
Moreover, we will denote by $\widetilde{\Xns/\iota_{\Xns}}$ (resp. $\widetilde{\Xns/\gamma}$) the minimal model of $\Xns/\iota_{\Xns}$ (resp. $\Xns/\gamma$).
\begin{proposition}
The automorphism $\iota_{\Xns}$ is the translation by the 2-torsion section on $\mathcal{E}:\Xns\ra\mathbb{P}^1$ and the automorphism $\gamma$ is the translation by a 4-torsion section. So:
\begin{itemize} \item 
	$k\left(\widetilde{\Xns/\iota_{\Xns}}\right)=k\left(\widetilde{\Xns/\gamma}\right)=k\left(\Xns\right)$;\item  $h^{i,0}\left(\widetilde{\Xns/\iota_{\Xns}}\right)=h^{i,0}\left(\widetilde{\Xns/\gamma}\right)=h^{i,0}\left(\Xns\right)$, for $i=0,1,2$; 
	\item  $\widetilde{\Xns/\iota_{\Xns}}$ admits an elliptic fibration $\mathcal{E}_{\iota}:\widetilde{\Xns/\iota_{\Xns}}\ra\mathbb{P}^1$ with Weierstrass equation
$y^2=x(x-1)(x-1+4\tau^n),$
three 2-torsion sections and singular fibers equal to $nI_2+I_{2n}+I_{2n}^*$ if $n$ is odd and to $nI_2+I_{2n}+I_{2n}$ if $n$ is even.
\item  $\widetilde{\Xns/\gamma}$ admits an elliptic fibration $\mathcal{E}_{\gamma}:\widetilde{\Xns/\gamma}\ra\mathbb{P}^1$ with Weierstrass equation
$y^2=x(x^2-2(4\tau^n+1)x+(4\tau^n-1)^2)$
one 4-torsion section, and singular fibers equal to $nI_4+I_n+I_n^*$ if $n$ is odd and to $nI_4+2I_{n}$ if $n$ is even. 
\end{itemize}
 
\end{proposition}
\proof
The automorphism $\iota_n\times\iota_n\in \Aut(\Cn\times \Cn)$ acts on the points fixed by $\alpha_n\times\alpha_n^{n-1}$ in the following way: $(1,1)\mapsto (-1,-1)$, $(1,-1)\mapsto (-1,1)$, $(\infty,\pm 1)\mapsto (\infty, \mp 1)$, $(\pm 1, \infty)\mapsto (\mp 1, \infty)$. Hence, the induced automorphism $\iota_{\Xns}$ preserves the classes $F_1$ and $F_2$ given in \eqref{eq: F1F2} and \eqref{eq: F2 even}. Moreover it acts on the curves which are supports of $F_1$ ($F_2$ respectively) as the translation by the 2-torsion section. So it is an involution which preserves two fibers and acts on each of them as a translation. There are two possibilities: either it acts trivially on the base of the fibration, and thus it is exactly the translation by a 2-torsion section, or it acts as an involution of the base of the fibration (which indeed preserves two points). In the latter case the involution does not preserve the 2-holomorphic forms, but $\iota_n\times\iota_n$ preserves all the 2-holomorphic form on $\Cn\times\Cn$, indeed $(\iota_n\times\iota_n)^*(\omega_j^{(1)}\wedge \omega_i^{(2)})=(-\omega_j^{(1)})\wedge(-\omega_i^{(2)})=\omega_j^{(1)}\wedge\omega_i^{(2)}$. We conclude that $\iota_{\Xns}$ is the translation by the 2-torsion section.
The quotient surface $\Xns/\iota_{\Xns}$ admits an elliptic fibration induced by the one of $\Xns$. So both the surfaces $\Xns$ and $\Xns/\iota_{\Xns}$ are birational to elliptic surfaces which are elliptic curves defined over the field of rational function $k(\tau)$ and the latter elliptic curve is obtained by the first as quotient by a translation by a 2-torsion rational point. This gives the Weierstrass equation of a birational model of $\Xns/\iota_{\Xns}$. The properties of the elliptic fibration are deduced from the Weiratrass equation.
The same arguments apply to the quotient $\Xns/\gamma$, once one proved that $\gamma$ is the translation by a 4-torsion section. First we observe that $\gamma$ preserves the 2-holomorphic form on $\Xns$, since  $(\varepsilon\circ(\iota_n\times\id))^*(\omega_j^{(1)}\wedge\omega_j^{(2)})=\varepsilon^*(-\omega_j^{(1)}\wedge\omega_j^{(2)})=\omega_j^{(1)}\wedge\omega_j^{(2)}$.
Moreover $\varepsilon\circ(\iota_n\times\id)$ acts on the points fixed by $\alpha_n\times\alpha_n^{n-1}$ in the following way: $(1,1)\mapsto (1,-1)\mapsto (-1,-1)\mapsto (-1, 1)$; $(\infty,1)\mapsto (1, \infty)\mapsto (\infty,-1)\mapsto (-1,\infty)$, and leaves invariant the point $(\infty,\infty)$. 
So $\gamma$ preserves the class $F_1$ and thus the elliptic fibration $|F_1|$ and acts on the reducible fibers $F_1$ and $F_2$ as a translation by a 4-torsion point. We observe that $\left(\varepsilon\circ(\iota_n\times\id)\right)^2=\iota_n\times\iota_n$. One concludes as above, recalling that the equation of the quotient by a 4-torsion point can be found in \cite[Proposition 2.2]{G}.

\subsection{The isotrivial fibrations on $\Xns$}
In Proposition \ref{prop: isotrivial firbation on Ynl} we observed that $\Yns$ admits two isotrivial fibrations whose generic fibers are isomorphic to $\Cn$. Since in the construction of $\Xns$ we do not contract any curve, each of these fibrations induces a fibration $\mathcal{F}_i:\Xns\ra\mathbb{P}^1$, $i=1,2$. If $n=3$ or $n=4$, the curve $\Cn$ is an elliptic curve (with complex multiplication) and the fibrations $\mathcal{F}_1$ and $\mathcal{F}_2$ are isotrivial elliptic fibrations. These fibrations do not coincide with the fibration $\mathcal{E}$, which in fact is not isotrivial. Indeed, if $n=3,4$, the K3 surfaces $\Xns$ admit more than one elliptic fibration. The K3 surfaces $X_{(3)}^{(2)}$ and $X_{(4)}^{(2)}$ are constructed in \cite{SI} exactly as quotient of the product of elliptic curves with complex multiplication and the presence of many elliptic fibrations on these K3 surfaces is well known. Indeed in \cite{Ni} the elliptic fibrations on both these K3 surfaces are classified. In particular one recognizes the fibration $\mathcal{E}$ and the isotrivial fibration $\mathcal{F}_i$ on $X_{(3)}^{(2)}$ (resp. $X_{(4)}^{(2)}$) as the fibrations number 6 and 5 in \cite[Table 1.1]{Ni} (resp. fibrations number 10 and 4 in  \cite[Table 1.2]{Ni}).

\section{The threefold $\Ynt$}\label{sec: 3-fold}
We observed in Proposition \ref{prop: Ynl quotients of products of Ynl} that in order to construct a variety which is birational to $\Ynt$ we can consider a birational model of  $\left(\Xns\times \Cn\right)/\left(\alpha_{\Xns}\times \alpha_n^{n-1}\right)$, where $\Xns$ is the smooth surface whose Weierstrass form is described in \eqref{cor: weierstrass Xns} and $\alpha_{\Xns}$ acts on $\tau$ as described in \eqref{eq: action alpha surface tau}. 
We will see that the threefold $\Ynt$ is obtained by a base change on a specific threefold which admits a fibration in K3 surfaces. This is in a certain sense the analogue of Proposition \ref{prop: Yns birat to base change}. The following theorem summarizes the main results of this section and it is a direct consequence of the Proposition \ref{prop: fibrations 3-folds n odd} and of the Corollary \ref{cor: Kodaira dim of 3-folds}.
\begin{theorem}\label{theorem: 3-folds}
The threefold $\Ynt$ is birational to a threefold $Z_{(n)}^{(3)}$ endowed with a fibration $\mathcal{Z}:Z_{(n)}^{(3)}\ra\mathbb{P}^1_t$. The minimal models of almost all the fibers are K3 surfaces $S_t$ with Picard number 19.

If $n=2$, $k\left(\Ynt\right)=-\infty$; if $n=3,4$, $\Ynt$ admits a desingularization which is a Calabi--Yau threefold (so $k\left(\Ynt\right)=0$); for any $n$, $k\left(\Ynt\right)\leq 1$.  
\end{theorem}

We first construct the fibration $\mathcal{Z}:Z_{(n)}^{(3)}\ra\mathbb{P}^1_t$ and we give its properties, in the following proposition. 
\begin{proposition}\label{prop: fibrations 3-folds n odd}
The threefold $\Ynt$ is birational to a threefold $Z_{(n)}^{(3)}$ endowed with a fibration $\mathcal{Z}:Z_{(n)}^{(3)}\ra\mathbb{P}^1_t$ whose generic fibers are birational to K3 surfaces $S_t$ with Picard number 19 and transcendental lattice $T_{S_t}\simeq U\oplus\langle 8\rangle$. There are $n$ special fibers (over $t$ such that $t^n=-1/4$) which are birational to K3 surfaces with Picard number 20 and transcendental lattice $T_{S_t}\simeq \langle 2\rangle\oplus \langle 8\rangle$ and two highly singular fibers (for $t=0,\infty$). The equation of the fibration $\mathcal{Z}$ is given by: 
\begin{equation}\label{eq: Weierstrass Ynt odd} Y^2=X\left(X^2+\left(v^2-1\right)\left(\left(v^2-1\right)-2 t^n  \right)X+t^{2n}\left(v^2-1\right)^2\right)\end{equation}
and $\mathcal{Z}$ is obtained by a base change of order $n$ from a fibration with equation: 
\begin{equation}\label{eq: Weierstrass Z1t} Y^2=X\left(X^2+\left(v^2-1\right)\left(\left(v^2-1\right)-2 t  \right)X+t^{2}\left(v^2-1\right)^2\right).\end{equation}
\end{proposition}
\proof
Let us consider the equation of the elliptic fibration $\Xns\ra\mathbb{P}^1_\tau$ obtained by putting $\sigma=1$ in \eqref{eq: weierstrass  Xns}, i.e. $y^2=x(x^2+(1-2\tau^n)x+\tau^{2n})$. The action of the automorphism $\alpha_{\Xns}$ is given in \eqref{eq: action alpha surface tau}. The functions $X:=xu^{2n}$ $Y:=yu^{3n}$ $t:=\tau u$ on $\Xns\times\Cn$ are invariant for the action of $\alpha_{\Xns}\times\alpha_n^{n-1}$. They satisfy the equation \eqref{eq: Weierstrass Ynt odd} and it is straightforward that there is a generically $n:1$ map between $\Xns\times\Cn$ and the threefold whose equation is \eqref{eq: Weierstrass Ynt odd}. Therefore, $\Ynt$ is birational to the threefold $Z_{(n)}^{(3)}$ defined in \eqref{eq: Weierstrass Ynt odd}. The equation \eqref{eq: Weierstrass Ynt odd} defines a fibration $\mathcal{Z}:Z_{(n)}^{(3)}\ra\mathbb{P}^1_{t}$ which is clearly induced by a base change of order $n$ (given by $t\mapsto t^n$ and thus branched at 0 and $\infty$) on the fibration \eqref{eq: Weierstrass Z1t}. Considering the curve $C_{(1)}$ defined as  $v^2=u+1$ (i.e. a rational curve double cover of $\mathbb{P}^1$), the threefold whose equation is \eqref{eq: Weierstrass Z1t} is obtained by the product of the rational elliptic surface $R$ (with equation $y^2=x(x^2+(1-2\tau)x+\tau^{2})$ and $C_{(1)}$, by using the 
functions $X:=xu^2$ $Y:=yu^3$ $t:=\tau u$.
Let us denote by $Z_{(1)}^{(3)}$ the threefold described by the equation \eqref{eq: Weierstrass Z1t}. The generic fiber $S_{\overline{t}}$ of the fibration $\mathcal{Z}:Z_{(1)}^{(3)}\ra\mathbb{P}^1_t$ over $\overline{t}$ is an  elliptic surface, whose Weierstrass equation is $$Y^2=X\left(X^2+\left(v^2-1\right)\left(\left(v^2-1\right)-2\overline{t}  \right)X+\overline{t}^{2}\left(v^2-1\right)^2\right),$$
considered as fibration to $\mathbb{P}^1_v$. The fundamental line bundle $\mathbb{L}$ of this Weierestrass equation has degree 2 (cf. \cite[Definitions II.4.1 and II.5.2]{M}). So the minimal model of the surface given by this Weierstrass equation is a K3 surface, by \cite[Lemma III.4.6]{M}. The discriminant of the Weierstrass equation $S_{\overline{t}}\ra\mathbb{P}^1$ is $$\Delta_{S_{\overline{t}}}=\overline{t}^{4}\left(v^2-1\right)^7\left(v^2-1-4\overline{t} \right).$$
If $\overline{t}\neq 0,\infty, -1/4$, then the singular fibers of the fibrations $S_{\overline{t}}\ra\mathbb{P}^1_v$ are two fibers of type $I_1^*$, over $v=\pm 1$, one fiber of type $I_8$ over $v=\infty$ and 2 fibers of type $I_1$ over $v=\pm (1+4\overline{t})$. Hence $\rho(S_{\overline{t}})\geq 19$. On the other hand the fibration $\mathcal{Z}:Z_{(1)}^{(3)}\ra\mathbb{P}^1_t$ is not an isotrivial fibration and $Z_{(1)}^{(3)}$ is at least a 1-dimensional family of K3 surfaces.  Indeed, if $\overline{t}=-\frac{1}{4}$, the reducible fibers of the fibration $S_{-1/4}\ra\mathbb{P}^1_v$ are $2I_1^*+I_8+I_2$ (and not $2I_1^*+I_8+2I_1$ as for the fibration on a general $S_{\overline{t}}$). We conclude that for a general $\overline{t}$, $\rho(S_{\overline{t}})=19$ and $\rho(S_{-1/4})=20$. In particular, this implies that, for a general $\overline{t}$, there are no sections of infinite order. The map $v\mapsto (X(v),Y(v))=(\overline{t}(v^2-1),\overline{t}(v^2-1)^2)$ is a  4-torsion section and the presence of a fiber $I_1^*$ prevents a higher order torsion section.
Hence $NS(S_{\overline{t}})$ is a specific overlattice of index 4 of $U\oplus D_5^{\oplus 2}\oplus A_7$ and the computation of its discriminant form implies that $T_{S_{\overline{t}}}\simeq U\oplus\langle 8\rangle$.  If $\overline{t}=-1/4$,  $NS(S_{-1/4})$ is a specific overlattice of index 4 of $U\oplus D_5^{\oplus 2}\oplus A_7\oplus A_1$ and the computation of its discriminant form implies that $T_{S_{-1/4}}\simeq \langle 2\rangle\oplus\langle 8\rangle$, see also \cite[Table 2]{SZ}. We observe that the fibers over $\overline{t}=0,\infty$ are no longer K3 surfaces.
The fibrations $Z_{(n)}^{(3)}\ra\mathbb{P}^1$ are obtained by a base change of order $n$ on the fibration $Z_{(1)}^{(3)}\ra\mathbb{P}^1$, so the results on the fibers of $Z_{(n)}^{(3)}\ra\mathbb{P}^1$ follows from the ones on the fibers of $Z_{(1)}^{(3)}\ra\mathbb{P}^1$.
\endproof
\begin{rem}{\rm
In the proof of Proposition \ref{prop: fibrations 3-folds n odd} we gave the explicit equation \eqref{eq: Weierstrass Ynt odd} for the threefold $Z_{(n)}^{(3)}$, which is a birational model of $\Ynt$. We observe that this equation can be read also in different way, indeed it is associated also to other two fibrations: \begin{itemize}
\item $Z_{(n)}^{(3)}\ra\mathbb{P}^1_{t}\times\mathbb{P}^1_{v}$, whose generic fibers are elliptic curves and whose discriminant is $\Delta=t^{4n}\left(v^2-1\right)^7\left(v^2-1-4t^{n}\right)$;
\item $Z_{(n)}^{(3)}\ra\mathbb{P}^1_{v}$ which is an isotrivial fibration whose generic fiber is birational to $\Xns$ (this is one of the fibrations described in Proposition \ref{prop: isotrivial firbation on Ynl}).
\end{itemize} }\end{rem}

In \cite{F}, the author studies a specific resolution $\Xnt$ of $\left(\Xns\times \Cn\right)/\left(\alpha_{\Xns}\times\alpha_n^{n-1}\right)$, in the case $n=3^c$ for some $c$. She proves that it is endowed with a fibration whose generic fiber, say $W$, is a K3 surface (and not only birational to a K3 surface). She also proves that the generic fiber has Picard number 19. The difference with our result is that we didn't fix a specific resolution and thus we can omit the condition $n=3^c$, but  we have to work up to birational equivalence.

In \cite{F} the author also considers certain specific curves on $W$ (rational or of genus 1), coming from the resolution she is considering. In order to compare the construction given in \cite[Section 6.2]{F} with the one introduced in the proof of the previous proposition, it could be interesting to have the full list of the elliptic fibrations on the K3 surface $S_{\overline{t}}$ (because this gives a large set of rational and genus 1 curves on $S_{\overline{t}}$). 
To obtain this classification it suffices to consider the method introduced in \cite{Ni} to the lattice $T_0\simeq A_7$ (where the lattice $T_0$ is chosen as required in \cite[Section 6.1]{Ni}).
The properties of the elliptic fibrations on $S_{\overline{t}}$ are listed in the following tables (the first line contains the reducible fibers, the second the rank of the Mordell--Weil group):
$$
\begin{array}{|c|c|c|c|c|c|c|c|c|c|c|c|}
\hline 
\mbox{red. fibers:}&2II^*&I_{12}^*&II^*+I_4^*&2III^*+2I_2&III^*+I_{10}&I_8^*+I_0^*&2I_4^*&I_{16}\\

\rk(MW):&1&1&1&1&1&1&1&2\\
\hline
\end{array}$$
$$
\begin{array}{|c|c|c|c|c|c|c|c|c|c|c|c|}
\hline 
\mbox{red. fibers:}&I_{5}^*+I_8&IV^*+I_3^*&I_2^*+I_{10}+I_2&2I_1^*+I_8&2I_9&I_{17}&I_{13}+I_{5}\\

\rk(MW):&	1&1&1&0&1&1&1\\
\hline
\end{array}
$$

In \cite{F}, the author finds two pairs of 9 rational curves on $W$, whose intersection is $\widetilde{A_8}$ and these curves could be the components of the fibration with reducible fibers $2I_9$ described in the table above.\\

The presence of the K3-fibration on $Z_{(n)}^{(3)}$, gives information on the Kodaira dimension of any birational model of this threefold. More precisely the following corollary holds.

\begin{corollary}\label{cor: Kodaira dim of 3-folds}
The Kodaira dimension of $\Ynt$ is at most 1. If $n=2$ the Kodaira dimension of $\Ynt$ is $-\infty$. If $n=3,4$, $\Ynt$ admits a desingularization which is a Calabi--Yau threefold (and in particular has Kodaira dimension equal to 0). 
\end{corollary}  
\proof We proved in Propositions \ref{prop: fibrations 3-folds n odd} that the threefold $Z_{(n)}^{(3)}$ admitting a fibration in K3 surfaces is birational to $\Ynt$. As in Proposition \ref{prop: Kodaira dimension}, the easy addiction formula implies  $$k(\Ynt)=k(Z_{(n)}^{(3)})\leq k(S_{\overline{t}})+\dim(\mathbb{P}^1)=1,$$ where $S_{\overline{t}}$ is a generic fiber of the fibration $\mathcal{Z}$ thus it is a K3 surface.

If $n=2$, the fibers of the isotrivial fibration $Z_{(n)}^{(3)}\ra \mathbb{P}^1_v$ are birational to the rational elliptic surface surface $X_{(2)}^{(2)}$ (see Proposition \ref{prop: the surface Z_n}). By applying the same formula as above to this fibration, one finds $k(Z_{(n)}^{(3)})\leq -\infty$. 

If $n=3$ or $n=4$ the $l$-fold $\Ynl$ is studied in \cite{CH}, where it is proved that it admits a crepant resolution which is a Calabi--Yau $l$-fold.\endproof

\section{The $l$-fold $\Ynl$ for $l>3$}\label{sec: l-fold}
We constructed the threefold $\Ynt$ from the surface $\Yns$ in the previous section. Iterating this process one is able to give an equation of a birational model of $\Ynl$ for any $l$. In particular this equation exhibits $\Ynl$ as birational to an elliptic fibration (see Corollary \ref{cor: elliptic fibration l-fold}) and to many other fibrations, some of which are isotrivial, given by projections on products of projective lines.
One of these fibrations is similar to the elliptic fibration on $\Xns$ and to the K3 fibration on a birational model of $\Ynt$, indeed it is a fibration in codimension 1 subvarieties and, even if we are not able to assure that the fibers are birational to Calabi--Yau $(l-1)$-folds, we can prove at least that the Kodaira dimension of the fiber is not bigger than 0. We use this property to bound the Kodaira dimension of $\Ynl$.

First we prove a useful result on the Kodaira dimension of the varieties described by certain Weierstrass equations.

\begin{lemma}\label{lemma: kod dim Weiestrass}
	Let $(v_i:w_i)$ be coordinates on the $i$-th copy of $\mathbb{P}^1$ in $(\mathbb{P}^1)^{\times l-1}$ and $a_i:=a_i((v_1:w_1),\ldots,(v_{l-1},w_{l-1}))$ be homogeneous polynomials of multidegree $(i,i,\ldots, i)$. 
	
	Let $W$ be an $l$-fold which admits an elliptic fibration $W\ra\left(\mathbb{P}^1_{v_i:w_1})\right)^{\times (l-1)}$ whose (extended) Weierstrass equation is \begin{equation}\label{eq: Weierstrass kodaira} y^2=x^3+a_{4}x^2+a_8x+a_{12}.\end{equation} Then, for any choice of the polynomials $a_i$, $k\left(W\right)\leq 0$.
\end{lemma}
\proof For a generic choice of the polynomials $a_i$, the equation \eqref{eq: Weierstrass kodaira} describes a smooth variety and the conditions on the multidegrees of the $a_i$'s implies that $W$ has a trivial canonical bundle (see \cite[Section (1.2)]{N} and observe that $\mathcal{L}^{-1}$ coincides with the canonical bundle of $\left(\mathbb{P}^1\right)^{\times(l-1)}$). In particular for generic choices of the polynomial $a_i$, $k(W)=0$. If the $a_i$'s are not general enough, $W$ is no longer necessarily smooth and the singularities may be not canonical, so a priori one has no control on the canonical bundle of the resolutions of $W$. In this case we consider three generic polynomials $b_4$, $b_8$ and $b_{12}$ in the same coordinates and with the same multidegree as the ones of $a_4$, $a_8$ and $a_{12}$ respectively. Let us consider the family $\mathcal{W}\ra \mathbb{P}^1_{\tau:\sigma}$ whose fiber over $(\tau:\sigma)$ is given by the equation
$$y^2=x^3+\left(\tau a_{4}+\sigma b_4\right)x^2+\left(\tau a_{8}+\sigma b_8\right)x+\left(\tau a_{12}+\sigma b_{12}\right).$$
The fibers over $(\sigma:\tau)=(0:1)$ is the variety $W$, and almost all the other fibers (in particular the ones over $(1:0)$) are smooth and thus their Kodiara dimension equals 0. Hence, by \cite[Theorem 1.2]{T}, one concludes that $k(W)\leq 0$.
\endproof
\begin{rem} {\rm The previous lemma (with the same proof) can be generalized taking as base $B$ of the elliptic fibration $W\ra B$ any product of projective spaces by requiring that the rational functions $a_i$ are sections of $|-2iK_B|$}.\end{rem}
\begin{theorem}\label{theorem: n-folds}
	The $l$-fold $\Ynl$ is birational to an $l$-fold $Z_{(n)}^{(l)}$ endowed with a fibration $\mathcal{Z}:Z_{(n)}^{(l)}\ra\mathbb{P}^1_t$ with the following equation
		\begin{equation}\label{eq: Weierstrass Ynl odd} Y^2=X\left(X^2+\prod_{i=1}^{l-2}\left(v_i^2-1\right)\left(\prod_{i=1}^{l-2}\left(v_i^2-1\right)-2 t^n  \right)X+t^{2n}\prod_{i=1}^{l-2}\left(v_i^2-1\right)^2\right)\end{equation}
	which is obtained by a base change of order $n$ from a fibration whose equation is: 
	$$Y^2=X\left(X^2+\prod_{i=1}^{l-2}\left(v_i^2-1\right)\left(\prod_{i=1}^{l-2}\left(v_i^2-1\right)-2 t\right)X+t^{2}\prod_{i=1}^{l-2}\left(v_i^2-1\right)^2\right).
	$$
	If $n=2$, $k\left(\Ynl\right)=-\infty$; if $n=3,4$, $\Ynl$ admits a desingularization which is a Calabi--Yau $l$-fold (so $k\left(\Ynl\right)=0$); for any $n$, $k\left(\Ynl\right)\leq 1$ if $l\geq 3$.  
	
\end{theorem}
\proof The statement on the presence of the fibration $\mathcal{Z}:Z_{(n)}^{(l)}\ra\mathbb{P}^1_t$ follows by iterated application of the method described in Proposition \ref{prop: fibrations 3-folds n odd}.
For a general $\overline{t}$, the fiber $F_{\overline{t}}$ of $\mathcal{Z}:Z_{(n)}^{(l)}\ra\mathbb{P}^1_t$ is 
$$Y^2=X\left(X^2+\prod_{i=1}^{l-2}\left(v_i^2-1\right)\left(\prod_{i=1}^{l-2}\left(v_i^2-1\right)-2 \overline{t}^n  \right)X+\overline{t}^{2n}\prod_{i=1}^{l-2}\left(v_i^2-1\right)^2\right).$$
This fiber is an $(l-1)$-fold edowed with an elliptic fibration $F_{\overline{t}}\ra \left(\mathbb{P}^1_{v_i}\right)^{\times (l-2)}$. The Weierstrass equation of this elliptic fibration satisfies the hypothesis of Lemma \ref{lemma: kod dim Weiestrass}, and thus $k\left(F_{\overline{t}}\right)\leq 0$. One concludes that $k\left(\Ynl\right)\leq 1$, as in the propositions \ref{prop: Kodaira dimension} and \ref{prop: fibrations 3-folds n odd}, by using the easy addiction formula.
The results for $n=3,4$ are proved in \cite{CH}. \endproof

\begin{corollary}\label{cor: elliptic fibration l-fold}
	The variety $\Znl$, which is birational to $\Ynl$, admits:\begin{itemize}\item an elliptic fibration $\mathcal{E}_Z:\Znl\ra \left(\mathbb{P}^1\right)^{\times l-1}$ given by the projection of \eqref{eq: Weierstrass Ynt odd}  to $\mathbb{P}^1_{t}\times \prod_{i=1}^{l-2} \mathbb{P}^1_{v_i}$; \item a K3-fibration given by the projection of \eqref{eq: Weierstrass Ynt odd}  to $\mathbb{P}^1_{t}\times \prod_{i=1}^{l-3} \mathbb{P}^1_{v_i}$;\item an isotrivial fibration with fibers birational to $Y_{(n)}^{(l-1)}$ given by the projection of \eqref{eq: Weierstrass Ynt odd} to $\mathbb{P}^1_{v_{l-2}}$.\end{itemize}  
	\end{corollary}

\end{document}